\documentclass[12pt]{article}
\usepackage{graphicx}  
\usepackage{amsmath}  
\usepackage{amssymb,amsthm}
\usepackage{color}
\usepackage{tikz}
\usepackage{pgfplots}

\newcommand{\R}{{\mathbb R}}

\newcommand{\Z}{{\mathbb Z}}\newcommand{\C}{{\mathbb C}}

\newcommand{\EE}{{\mathbf E}}
\newcommand{\FF}{{\mathbf F}}
\newcommand{\ee}{{e}}
\newcommand{\ii}{{i}}

\let\epsilon\varepsilon
\let\theta\vartheta
\let\hat\widehat

\newtheorem{theorem}{Theorem}[section]\newtheorem{lemma}[theorem]{Lemma}

\newtheorem{corollary}[theorem]{Corollary}

\newtheorem{remark}[theorem]{Remark}

	\title{NLS approximation for a scalar FPUT system 
		on a 2D square lattice with a cubic nonlinearity}
	
	\author{Ioannis Giannoulis$^{1} $,   Bernd Schmidt$^{2} $   and Guido Schneider$^{3} $\\
		{\small
			$^{1}$ Department of Mathematics,
University of Ioannina,} 
{\small University Campus, 45110 Ioannina, Greece}
\\
{\small $^{2}$  Institut f\"ur Mathematik, Universit\"at  Augsburg, 	Universit\"atsstra\ss{}e 14, 86159 Augsburg, Germany
}\\
		{\small
			$^{3}$ Institut f\"ur Analysis, Dynamik und Modellierung, Universit\"at Stuttgart, 70569 Stuttgart, Germany}}

\begin{document}
	
\maketitle
	
	\begin{abstract}
		We consider a scalar Fermi-Pasta-Ulam-Tsingou (FPUT) system on a square 2D lattice with a cubic nonlinearity.  For such systems the  NLS equation   can be derived to describe 
		the evolution  of an oscillating moving wave packet 		
		 of small amplitude which is  slowly modulated  in time and space. We show that this NLS approximation makes correct predictions about the dynamics of the original scalar FPUT system for the strain and the displacement variables. 
	\end{abstract}

\section{Introduction}

The equations of motion of 
a scalar Fermi-Pasta-Ulam-Tsingou (FPUT) system on a two-dimensional  square  lattice
are given by 
\begin{eqnarray}  
\label{fpuintro}
\partial_t^2 q_{m,n} & = & W'(q_{m+1,n}-q_{m,n})- W'(q_{m,n}-q_{m-1,n})\\ && +W'(q_{m,n+1}-q_{m,n})- W'(q_{m,n}-q_{m,n-1})
\nonumber
\end{eqnarray}
for all $ (m,n) \in \Z^2 $.  
The variables $ q_{m,n} \in \R $ describe a vertical displacement 
in $ z $-direction of a particle with unit mass located at the $(m,n)$-th site 
of the lattice in the $(x,y)$-plane.  In this paper we consider 
\eqref{fpuintro} with an interaction force $ W'$
with linear, cubic and higher order  terms but no quadratic terms. W.l.o.g. 
for our purposes in this paper we assume
\begin{equation} \label{Wprime}
W'(u) = u -u^3.
\end{equation}
For \eqref{fpuintro} an NLS equation  can be derived in order to describe 
		the evolution  of an oscillating moving wave packet 		
		 of small amplitude which is  slowly modulated  in time and space
		 by some envelope. 
We are interested in the validity of  this approximation. 

We can derive the NLS equation for the original displacement variables 
$ q_{m,n} $ or the associated strain variables. 
In this paper we start with the strain variables and come back to the 
handling of the displacement variables in Section \ref{secdisc}.
We follow \cite{PS23}, where the KP approximation has been justified for 
a quadratic FPUT system,
and introduce 
the strain variables 
\begin{equation} \label{strain}
u_{m,n} = q_{m+1,n}-q_{m,n}, \qquad v_{m,n} = q_{m,n+1}-q_{m,n}
\end{equation}
which leads to the compatibility condition 
\begin{equation} \label{compiintro}
u_{m,n+1}-u_{m,n} = v_{m+1,n}-v_{m,n} .
\end{equation}
The  strain variables satisfy
\begin{eqnarray} 
 \label{umneq} 
\partial_t^2 u_{m,n} & = &  W'(u_{m+1,n})- 2 W'(u_{m,n}) + W'(u_{m-1,n})
\\ &&
\qquad + W'(v_{m+1,n})- W'(v_{m+1,n-1}) -W'(v_{m,n})+ W'(v_{m,n-1}), 
\nonumber
\\
\label{vmneq} 
\partial_t^2 v_{m,n} & = &  W'(v_{m,n+1})- 2 W'(v_{m,n}) + W'(v_{m,n-1})  
\\ && \qquad + W'(u_{m,n+1})- W'(u_{m-1,n+1}) -W'(u_{m,n})+ W'(u_{m-1,n}).
\nonumber
\end{eqnarray}
By undoing the introduction of $  u_{m,n} $ and $ v_{m,n} $ through \eqref{strain},
this system is equivalent to the original FPUT system \eqref{fpuintro}
for the displacements. Moreover, the compatibility condition \eqref{compiintro}
 is invariant with respect to the time evolution of \eqref{umneq}-\eqref{vmneq}, 
see \eqref{uvrelation} and cf. \cite{PS23}. 

In this introduction we summarize the derivation of the NLS equation and state our 
result. Details can be found in subsequent sections.
The NLS ansatz for the strain  variables is given by  
\begin{eqnarray}
\psi_{u,m,n}(t) & = & \varepsilon A(X,Y,T) e^{i (k_0 m + l_0 n + \omega_0 t) } + c.c., \label{uapprox1} \\ 
\psi_{v,m,n}(t) & = & \varepsilon  B(X,Y,T) e^{i (k_0 m + l_0 n + \omega_0 t) } + c.c. , \label{vapprox1}
\end{eqnarray} 
with
\begin{equation} \label{XYT}
X =  \varepsilon(m + c_x t),  
\quad Y = \varepsilon( n+c_y t),  \quad T = \varepsilon^2 t, \quad 0 < \varepsilon \ll 1 
\end{equation}
where in lowest order $ A $ and $ B $ are related through an expansion of the 
compatibility condition which leads to
\begin{equation} \label{ABrelation}
(e^{ik_0}-1)  B(X,Y,T) =  (e^{il_0}-1) A(X,Y,T) , 
\quad t \geq 0.
\end{equation}
Inserting the ansatz \eqref{uapprox1}-\eqref{vapprox1} into the system \eqref{umneq}-\eqref{vmneq} and 
equating the coefficients in front of the $   \varepsilon^j e^{i (k_0 m+l_0n  +  \omega_0 t) } $ to zero 
gives  the dispersion relation $ \omega_0  = \omega  (k_0,l_0) $ with $ \omega $ given in \eqref{omega} 
for $ j =1 $,
 the identity $ (c_x,c_y) = \nabla \omega (k_0, l_0) $ for the group velocity
for $ j = 2 $, 
and finally for $ j = 3 $ the NLS equation
\begin{equation} \label{NLSintro}
\partial_T A  =  - \frac12 \ii (\partial_X,\partial_Y) \nabla^2 \omega(\vec{k}_0) (\partial_X,\partial_Y)^T A 
+ 
4 
\gamma_A  |A|^2 A  
\end{equation}
with 
$\gamma_A$  as in 
\eqref{gammaA} below. 
For the local existence and uniqueness theory for \eqref{NLSintro} we refer to \cite{Cazenave}. 
The goal of this paper is to prove the following approximation theorem.
\begin{theorem} \label{th11}
Let $ s_A > 4 $ and $ (k_0,l_0) \neq (0,0) $ chosen in such a way that the subsequent 
non-resonance condition \eqref{nonres} is satisfied.
For all $ T_0 > 0 $, $ C_1 > 0 $, $ C_2 > 0 $  
there exist $ \varepsilon_0 > 0 $ and $ C_3 > 0 $ such that for all 
$ \varepsilon \in (0,\varepsilon_0) $ the following holds.
In case $ k_0 \neq 0 $
let $ A \in C([0,T_0],H^{s_A}(\R^2,\C)) $ be a solution of the  NLS equation \eqref{NLSintro} with 
$$ 
\sup_{T \in [0,T_0]} \| A(\cdot,T) \|_{H^{s_A}} \leq C_1 
$$ 
and let $\psi_{u,m,n}$, $\psi_{v,m,n}$ be defined by \eqref{uapprox1}-\eqref{vapprox1} with $B$ given by \eqref{ABrelation}. 
In case $ k_0 = 0 $ consider $ B $ instead of $ A $. 
Then for
all initial conditions of \eqref{umneq}-\eqref{vmneq} 
satisfying the  compatibility condition \eqref{compiintro} and
\begin{eqnarray*}
\lefteqn{ \sup_{(m,n) \in \Z^2} (| u_{m,n}(0) -  \psi_{u,m,n}
(0) | +| v_{m,n}(0) -  \psi_{v,m,n}
(0) |  } \\ && +| \partial_t u_{m,n}(0) -  \partial_t \psi_{u,m,n}
(0) |+ | \partial_t v_{m,n}(0) -  \partial_t \psi_{v,m,n}
(0) |)  \leq  C_2 \varepsilon^2 
\end{eqnarray*} 
the
solutions $ (u_{m,n},v_{m,n})_{(m,n) \in \Z^2} $ of 
\eqref{umneq}-\eqref{vmneq} with these initial conditions  satisfy
\begin{eqnarray*} 
\lefteqn{ \sup_{t \in [0,T_0/\varepsilon^2]} \sup_{(m,n) \in \Z^2} (| u_{m,n}(t) -  \psi_{u,m,n}
(t) | +| v_{m,n}(t) -  \psi_{v,m,n}
(t) |}
\\ && +| \partial_t u_{m,n}(t) -  \partial_t \psi_{u,m,n}
(t) |+ | \partial_t v_{m,n}(t) -  \partial_t \psi_{v,m,n}
(t) |
 ) \leq C_3 \varepsilon^2.
 \end{eqnarray*}  
\end{theorem}
\begin{remark}{\rm
As we will see below it is always possible to satisfy the assumptions of Theorem \ref{th11}, i.e., to find initial conditions which are $ \mathcal{O}(\varepsilon^2) $-close to the approximation and satisfy the compatibility condition.
}\end{remark}

\begin{remark}{\rm
A corresponding approximation result for the original displacement variables $ q_{m,n} $ can be obtained in an analogous way, cf.\ Theorem~\ref{th51} below. 
}\end{remark}

Over the last decades the derivation of amplitude equations and their justification for one-dimensional atomic chains has received a lot of attention. In particular, for mono- and poly-atomic FPUT models KdV approximation results have been obtained for instance in \cite{SW00fpu,CCPS12,GMWZ14,HKY21} and NLS approximation results in \cite{GM04,GM06,Schn10}. In contrast, corresponding approximation results for amplitude equations on lattices beyond the one-dimensional set-up are considerably less well explored. Among the exceptions are the two recent contributions \cite{HP22,PS23}, in which KP approximation results have been established for FPUT models in 2D, and also \cite{GHM06}, in which an 
NLS approximation result is shown for higher dimensional lattices, albeit with the additional assumption of an on-site potential that 
anchors individual atoms in the background. 

To the best of our knowledge, the present paper provides the first NLS approximation result for a 2D lattice in a full Galilean invariant setting without imposing an on-site potential. Moreover both the original displacement and the strain variables are considered. In the following, we restrict our analysis to cubic nonlinearities. From a technical point of view, this simplifies the final error estimates since a simple application of Gronwall's inequality is sufficient to obtain the error estimates on the long $ \mathcal{O}(1/\varepsilon^2) $-time scale. More importantly, this set-up will enable us in the future to contrast the presented approximation results for the displacement and strain variables with the general case of non-vanishing quadratic nonlinearities. In this case, the Davey-Stewartson system will take over the role of the NLS equation. The analysis of this subtle regime is the subject of ongoing work. 

The paper is organized as follows. It turns out to be advantageous 
to work in Fourier space. Therefore, in 
the next section we derive the Fourier transformed version of \eqref{umneq}-\eqref{vmneq}
and write it  as a 
first order system. 
In Section \ref{sec3} we derive the NLS equation 
and estimate the residual terms. In Section \ref{sec4} we use 
Gronwall's inequality to estimate the error made by this approximation.
The paper is closed with a discussion section where we 
discuss an approximation result for the original displacement variables $ q_{m,n} $ 
and the situation of a small uncertainty in the interaction forces.
\medskip

{\bf Notation.} 
Possibly different constants are denoted with the same symbol $ C $ if they can be chosen independent of the small perturbation parameter $ 0 < \varepsilon^2 \ll 1 $.
We define $ \|u \|_{H^s(\R^2)} = \|\widehat{u} \|_{L^2_s(\R^2)} $, 
where the Fourier transform $ \widehat{u} $ of $ u $ is defined in \eqref{FT} 
and where for $ s \ge 0 $ and $p \in [1, \infty)$  the space $ L^p_s $ 
is equipped with the norm $ \| \widehat u \|_{L^p_s(\R^2)} = \| \widehat u \rho^s \|_{L^p(\R^2)} $ 
with $ \rho(k,l) = (1+k^2+l^2)^{1/2} $. 
\medskip

{\bf Acknowledgement.}  
The work of  Ioannis Giannoulis was supported by the Bilateral Exchange of Academics 
2021 and 2024 Programs   
between Germany (DAAD) and Greece (IKY).
The work of Guido Schneider  is supported by the Deutsche Forschungsgemeinschaft
DFG through the cluster of excellence  'SimTech'  under EXC 2075-390740016.

\section{The Fourier transformed FPUT system}

In order to specify the compatibility condition and to derive the NLS equation 
we will work in Fourier space, 
like for the FPUT system on the 1D lattice \cite{Schn10}. 
 Therefore, we define
$$
\hat{u}(k,l) = \frac{1}{(2\pi)^2} \sum_{(m,n) \in \mathbb{Z}^2} u_{m,n} e^{-ikm-iln}, \qquad 
u_{m,n} = 
\int_{\mathbb{T}^2} \hat{u}(k,l) e^{ikm+iln}  d(k,l), 
$$ 
and similarly for $v_{m,n}$ and $w_{m,n} = \partial_t q_{m,n}$, where $\mathbb{T} := (-\pi,\pi]$ is equipped with periodic boundary conditions. 
From the definition \eqref{strain} of the strain variables we obtain in Fourier space
$$
\partial_t \hat{u}(k,l,t) = (e^{ik} - 1) \hat{w}(k,l,t), \quad \partial_t \hat{v}(k,l,t) = (e^{il} - 1) \hat{w}(k,l,t),
$$
which implies that the compatibility condition 
\begin{equation} \label{uvrelation}
(e^{ik}-1)  \widehat{v}(k,l,t)  =  (e^{il}-1) \widehat{u}(k,l,t), 
\quad t \geq 0, 
\end{equation}
is invariant with 
respect to the time evolution of the FPUT system \eqref{umneq}-\eqref{vmneq}.
The linearized system 
\begin{eqnarray*} 
\partial_t^2 u_{m,n} & = & u_{m+1,n}- 2 u_{m,n} + u_{m-1,n} + v_{m+1,n}- v_{m+1,n-1} -v_{m,n}+ v_{m,n-1}, \\
\partial_t^2 v_{m,n} & = & v_{m,n+1}- 2 v_{m,n} + v_{m,n-1}
+ u_{m,n+1}- u_{m-1,n+1} -u_{m,n}+ u_{m-1,n},
\end{eqnarray*}
is written in Fourier space as
\begin{eqnarray*} 
\partial_t^2 \widehat{u}(k,l,t)  & = &  (e^{i k}-2+e^{-i k}) \widehat{u}(k,l,t)  + (e^{ik}-1)(1-e^{-il})\widehat{v}(k,l,t), \\
\partial_t^2 \widehat{v}(k,l,t) & = & (e^{i l}-2+e^{-i l}) \widehat{v}(k,l,t)  + (e^{il}-1)(1-e^{-ik})\widehat{u}(k,l,t), 
\end{eqnarray*}
where we have used
$$ 
e^{ik}- e^{ik}e^{-il}-1 + e^{-il} = e^{ik}(1-e^{-il}) - (1-e^{-il})= (e^{ik}-1)(1-e^{-il}).
$$ 
With 
\begin{equation} \label{omega}
\omega_x^2(k) := 2-e^{-i k}-e^{i k}, \qquad 
\omega_y^2(l) := 2-e^{-i l}-e^{i l}, \qquad
\omega(k,l) { : = } (\omega_x(k)^2 +  \omega_y(l)^2)^{1/2} ,
\end{equation}
and the compatibility condition \eqref{uvrelation},
 the linearized system can be written as
\begin{equation}
\label{linear-pr}
\partial_t^2 \widehat{u} + \omega^2 \widehat{u} = 0 \quad 
\textrm{and} \quad 
\partial_t^2 \widehat{v} + \omega^2  \widehat{v} = 0 .
\end{equation} 
Extending exactly the same calculations  
for $W'(u) = u - u^3$ instead of $W'(u) = u $, 
we obtain that the nonlinear FPUT system  \eqref{umneq}-\eqref{vmneq} takes in Fourier space the form
\begin{eqnarray*} 
\partial_t^2 \widehat{u}(k,l,t) & = &  -\omega_x^2(k)(\widehat{u} - \widehat{u}*\widehat{u}*\widehat{u})(k,l,t)+ (e^{ik}-1)(1-e^{-il})(\widehat{v}- \widehat{v}*\widehat{v}*\widehat{v})(k,l,t), \\
\partial_t^2 \widehat{v}(k,l,t) & = &  -\omega_y^2(l)  (\widehat{v} - \widehat{v}*\widehat{v}*\widehat{v})(k,l,t) + (e^{il}-1)(1-e^{-ik})(\widehat{u} - \widehat{u}*\widehat{u}*\widehat{u})(k,l,t)
\end{eqnarray*}
where 
\begin{equation} \label{perconv}
\widehat{(u v)}(k,l) = 
(\widehat{u}*\widehat{v})(k,l) =
 \int_{\mathbb{T}^2}   
 \widehat{u}(k-\widetilde{k},l- \widetilde{l})\widehat{v}(\widetilde{k},\widetilde{l}) 
d ( \widetilde{k} , \widetilde{l} ) .
\end{equation}
By using the compatibility condition \eqref{uvrelation}
and the notations $$ \rho_u(k,l) = (e^{ik}-1)(1-e^{-il}) \quad \textrm{and}  \quad 
 \rho_v(k,l) =(e^{il}-1)(1-e^{-ik}) $$ 
we rewrite this system in the form:
\begin{eqnarray}  \label{syst1}
\partial_t^2 \widehat{u} & = & -\omega^2  \widehat{u} + \omega_x^2
(\widehat{u}*\widehat{u}*\widehat{u}) - \rho_u (\widehat{v}*\widehat{v}*\widehat{v}), 
\\  
\partial_t^2 \widehat{v} & = & -\omega^2  \widehat{v} + \omega_y^2(\widehat{v}*\widehat{v}*\widehat{v})- \rho_v ( \widehat{u}*\widehat{u}*\widehat{u}). \label{syst2}
\end{eqnarray}
This system in combination with the compatibility condition \eqref{uvrelation} is the starting point for the derivation and justification 
of the NLS equation.  We write  \eqref{syst1}-\eqref{syst2} as  a first order system
\begin{eqnarray*} 
\partial_t \widehat{u} & = & i \omega \widehat{u}_1 , \\
\partial_t \widehat{u}_1 & = & i \omega \widehat{u} + \frac{\omega_x^2}{i \omega}
(\widehat{u}*\widehat{u}*\widehat{u}) - \frac{\rho_u}{i \omega}
(\widehat{v}*\widehat{v}*\widehat{v}), 
\\  
\partial_t \widehat{v} & = & i \omega \widehat{v}_1 , \\
\partial_t \widehat{v}_1 & = & i \omega \widehat{v} 
+ \frac{\omega_y^2}{i \omega}(\widehat{v}*\widehat{v}*\widehat{v})- \frac{\rho_v}{i \omega}( \widehat{u}*\widehat{u}*\widehat{u}),
\end{eqnarray*}
with 
$$
\frac{\omega_x^2(k)}{i \omega(k,l)} = \mathcal{O}(\frac{k^2}{(k^2+l^2)^{1/2}}), \qquad 
 \frac{\rho_u(k,l)}{i \omega(k,l)} = \frac{(e^{ik}-1)(1-e^{-il})}{i \omega(k,l)} = \mathcal{O}(\frac{kl}{(k^2+l^2)^{1/2}}), \qquad \textrm{etc.}
$$
which are bounded for $ (k,l) \to (0,0) $.
We diagonalize the last system by 
$ \widehat{U}_{1} = \widehat{u} + \widehat{u}_{1} $, $ \widehat{U}_{-1} = \widehat{u} - \widehat{u}_{1} $
and similar 
for the $ v $-variables. We find 
\begin{eqnarray} \label{syst4}
\partial_t \widehat{U}_1 & = & i \omega \widehat{U}_1 + \frac{\omega_x^2}{8 i \omega}
(\widehat{U}_1 + \widehat{U}_{-1})^{*3} - \frac{\rho_u}{8 i \omega}
(\widehat{V}_1 + \widehat{V}_{-1})^{*3}
, \\
\partial_t \widehat{U}_{-1} & = & - i \omega \widehat{U}_{-1} - \frac{\omega_x^2}{8 i \omega}
(\widehat{U}_1 + \widehat{U}_{-1})^{*3} + \frac{\rho_u}{8 i \omega}
(\widehat{V}_1 + \widehat{V}_{-1})^{*3}, 
\\  
\partial_t \widehat{V}_1 & = & i \omega \widehat{V}_1 + \frac{\omega_y^2}{8 i \omega}(\widehat{V}_1 + \widehat{V}_{-1})^{*3}- \frac{\rho_v}{8 i \omega}(\widehat{U}_1 + \widehat{U}_{-1})^{*3} , \\
\partial_t \widehat{V}_{-1} & = & - i \omega \widehat{V}_{-1} 
- \frac{\omega_y^2}{8 i \omega}(\widehat{V}_1 + \widehat{V}_{-1})^{*3}+ \frac{\rho_v}{8 i \omega}(\widehat{U}_1 + \widehat{U}_{-1})^{*3} ,
\label{syst5}
\end{eqnarray}
where $ \widehat{U}^{*3} = \widehat{U}*\widehat{U}*\widehat{U} $.
In contrast to the second order system the coefficients in the first order 
system are not smooth at the wave vector $ (k,l) = (0,0) $.
However, since we have a cubic nonlinearity,  subsequently for the NLS approximation this wave vector will not play a role. The compatibility condition for \eqref{syst4}-\eqref{syst5} is obtained 
as follows. From \eqref{uvrelation} and the definition of $ u_1 $ and $ v_1 $ 
we obtain that 
\begin{equation} \label{u1v1relation}
(e^{ik}-1)  \widehat{v}_1(k,l,t)  =  (e^{il}-1) \widehat{u}_1(k,l,t), 
\quad t \geq 0.
\end{equation}
After the diagonalization we obtain 
\begin{equation} \label{UVrelation}
(e^{ik}-1)  \widehat{V}_j(k,l,t)  =  (e^{il}-1) \widehat{U}_j(k,l,t), 
\quad t \geq 0, \quad j \in \{-1,1\}.
\end{equation}

\section{Derivation of the NLS equation}
\label{sec3}

i) A serious difficulty which occurs is the fact that the solutions of the FPUT
problem \eqref{umneq}-\eqref{vmneq} live on $ \Z^2 $ but the 
solutions of the NLS equation \eqref{NLSintro} live on $ \R^2 $.
Hence we also need the Fourier  transform on $ \R^2 $ and an operator 
which connects both.

For a function $ A:\R^2 \to \C$ we denote by $ \widehat{A} $ its Fourier transform, i.e., 
\begin{equation} \label{FT}
\hat{A}(K,L) = \frac{1}{(2\pi)^2} \int_{\mathbb{R}^2} A(X,Y) e^{-i KX-iLY} d(X,Y), 
\quad 
A(X,Y) = \int_{\mathbb{R}^2} \hat{A}(K,L) e^{i K X + i L Y} d(K, L).
\end{equation} 
With $ \widehat{A}^{\varepsilon} $ we denote the restriction of $ \widehat{A} $ 
to $ (-\pi/\varepsilon,\pi/\varepsilon]^2  $ and by 
$ \widehat{A}^{\varepsilon,per}$  the $ 2 \pi/\varepsilon $-periodic continuation 
of $ \widehat{A}^{\varepsilon} $ which lives on $ \mathbb{T}^2_{2 \pi/\varepsilon} $, 
where $ \mathbb{T}_{2 \pi/\varepsilon}  := (-\pi/\varepsilon,\pi/\varepsilon]$.
%

ii) In this section we derive the NLS equation 
for the first order system  \eqref{syst4}-\eqref{syst5}. 
We extend the ansatz for the derivation of the NLS equation 
in Fourier space that corresponds to \eqref{uapprox1}-\eqref{vapprox1} in physical space 
by higher order terms in order to make the residual sufficiently small.
The residual contains all terms which do not cancel after inserting the ansatz 
into the equations. In order to measure the size of the residual  in Fourier space we use the space 
$ L^1(\mathbb{T}^2) $. This space is closed under convolution, and 
multiplication operators $ \widehat{\Lambda} =   \widehat{\Lambda}(k) $
can be estimated in this space by 
$$ 
\|\widehat{\Lambda}\widehat{u} \|_{L^1} \leq\|\widehat{\Lambda} \|_{L^{\infty}} \| \widehat{u} \|_{L^1} .
$$  
iii) In the following we use the notation $\vec{k} = (k,l)$. 
First we  choose a wave vector  
$ \vec{k}_0 = (k_0,l_0) \in (\mathbb{T} \setminus \{ 0 \} )^2$ 
and then come back to $ \vec{k}_0 = (k_0,l_0) \in \mathbb{T}^2 \setminus \{ \vec{0} \} $
in Remark \ref{remAB}. 
Since the Fourier transform of $ \varepsilon A(\varepsilon x, \varepsilon y) $  
w.r.t.\ $ x , y \in \R $ is given by 
$ \varepsilon^{-1} \hat{A}(\frac{\vec{k}}{\varepsilon}) $
%
we make the ansatz for  $\widehat{U}_1 $, $ \widehat{U}_{-1} $:
\begin{eqnarray} \label{eq21a}
\hat{\psi}_{U,1}(\vec{k},t) & = & \varepsilon^{-1} 
\hat{A}^{\varepsilon,per}_1(\frac{\vec{k}-\vec{k}_0}{\varepsilon},\varepsilon^2 t)
\EE \FF_1 +  \varepsilon \hat{A}^{\varepsilon,per}_{1,-1}(\frac{\vec{k}+\vec{k}_0}{\varepsilon},\varepsilon^2 t)
\EE^{-1}  \FF_{-1}  \\ && 
+ \varepsilon \hat{A}^{\varepsilon,per}_{1,3}(\frac{\vec{k}-3 \vec{k}_0}{\varepsilon},\varepsilon^2 t)
\EE^3 \FF_3 
+ \varepsilon \hat{A}^{\varepsilon,per}_{1,-3}(\frac{\vec{k}+3 \vec{k}_0}{\varepsilon},\varepsilon^2 t)
\EE^{-3} \FF_{-3},  \nonumber 
\\
\label{eq21b}
\hat{\psi}_{U,-1}(\vec{k},t) & = &
 \varepsilon^{-1} \hat{A}^{\varepsilon,per}_{-1}(\frac{\vec{k}+\vec{k}_0}{\varepsilon},\varepsilon^2 t)
\EE^{-1} \FF_{-1} 
+ \varepsilon \hat{A}^{\varepsilon,per}_{-1,1}(\frac{\vec{k}-\vec{k}_0}{\varepsilon},\varepsilon^2 t)
\EE \FF_1  \\ && 
+ \varepsilon \hat{A}^{\varepsilon,per}_{-1,3}(\frac{\vec{k}-3 \vec{k}_0}{\varepsilon},\varepsilon^2 t)
\EE^{3} \FF_3 
+ \varepsilon \hat{A}^{\varepsilon,per}_{-1,-3}(\frac{\vec{k}+3 \vec{k}_0}{\varepsilon},\varepsilon^2 t)
\EE^{-3} \FF_{-3},  
\nonumber 
\end{eqnarray}
where
$$ \EE =
\ee^{\ii \omega(\vec{k}_0)t} \quad \textrm{ and } \qquad   \FF_m = \ee^{\ii \nabla\omega(\vec{k}_0)\cdot (\vec{k}-m \vec{k}_0) t}      \quad \textrm{for }  m \in \{-3,-1,1,3\} ,$$
with $ \nabla\omega(\vec{k}_0) = (c_x,c_y) $ in \eqref{XYT} and $ A_{-1} = \overline{A_1} $, $ A_{-1,1} = \overline{A_{1,-1}} $, $ A_{-1,3} = \overline{A_{1,-3}} $,
and $ A_{-1,-3} = \overline{A_{1,3}} $
in physical space, 
and similar for $ \hat{\psi}_{V,1} $, $ \hat{\psi}_{V,-1} $,  where in their definition 
$ A $ is replaced by $ B $.

iv) Inserting this ansatz into \eqref{syst4}-\eqref{syst5}  gives the 
approximation equations which we compute by using the following facts.
 The nonlinear terms below can be computed by using 
$$
(\widehat{U}_1 + \widehat{U}_{-1})^{*3} = \widehat{U}_1^{*3} 
+ 3  \widehat{U}_1^{*2} *  \widehat{U}_{-1}
+ 3  \widehat{U}_1*  \widehat{U}_{-1}^{*2}  +
\widehat{U}_{-1}^{*3}.
$$
Rescaling the convolution  integrals 
and taking formally the limit 
$\int_{\mathbb{T}^2_{2\pi/\varepsilon}} \int_{\mathbb{T}^2_{2\pi/\varepsilon}}\to \int_{\R^2} \int_{\R^2}$ 
gives for instance
\begin{eqnarray*}
\lefteqn{
( \frac{3 \omega_x^2}{8 i \omega}
\widehat{U}_1 * \widehat{U}_1 * \widehat{U}_{-1})(\vec{k},t) 
 =  \frac{3 {  \omega_x^2(k) }}{8 i \omega(\vec{k})} 
\int_{\mathbb{T}^2} \int_{\mathbb{T}^2}  \widehat{U}_1(\vec{k}-\vec{k}_1,t)  \widehat{U}_1(\vec{k}_1-\vec{k}_2,t)  \widehat{U}_{-1}(\vec{k}_2,t) d \vec{k}_2 d \vec{k}_1} 
\\ & = & \frac{3  \omega_x^2(k)  }{8 i \omega(\vec{k})} 
\int_{\mathbb{T}^2} \int_{\mathbb{T}^2} \varepsilon^{-1} \hat{A}^{\varepsilon,per}_{1}(\frac{\vec{k}-\vec{k}_0-\vec{k}_1}{\varepsilon},T)
 \varepsilon^{-1} \hat{A}^{\varepsilon,per}_{1}(\frac{\vec{k}_1-\vec{k}_0-\vec{k}_2}{\varepsilon},T) \\ && \qquad \qquad \times
\varepsilon^{-1} \hat{A}^{\varepsilon,per}_{-1}(\frac{\vec{k}_2+\vec{k}_0}{\varepsilon},T)
d \vec{k}_2 d \vec{k}_1 { \EE \FF_1 } + h.o.t.
\\& = &\varepsilon  \frac{3 \omega_x^2(k_0+ \varepsilon K )  }{8 i \omega(\vec{k_0} + \varepsilon \vec{K})} 
\int_{\mathbb{T}^2_{2\pi/\varepsilon}} \int_{\mathbb{T}^2_{2\pi/\varepsilon}}  \hat{A}^{\varepsilon,per}_{1}(\vec{K}-\vec{K}_1,T)
 \hat{A}^{\varepsilon,per}_{1}(\vec{K}_1-\vec{K}_2,T) \\ && \qquad \qquad \times
\hat{A}^{\varepsilon,per}_{-1}(\vec{K}_2,T)
d \vec{K}_2 d \vec{K}_1 {  \EE  \ee^{\ii \nabla\omega(\vec{k}_0)\cdot \varepsilon \vec{K} t}  }  + h.o.t.
\\ &\to &
\varepsilon  \frac{3   \omega_x^2(k_0)  }{8 i \omega(\vec{k_0} )} 
\int_{\R^2} \int_{\R^2} \hat{A}_{1}(\vec{K}-\vec{K}_1,T)
 \hat{A}_{1}(\vec{K}_1-\vec{K}_2,T) 
\hat{A}_{-1}(\vec{K}_2,T)
d \vec{K}_2 d \vec{K}_1 { \EE } + h.o.t.
\end{eqnarray*}
in the limit $ \varepsilon \to 0 $.
Furthermore, using 
$$
\partial_t ( \varepsilon^{-1} \hat{A}^{\varepsilon,per}_1(\frac{\vec{k}-\vec{k}_0}{\varepsilon},\varepsilon^2 t) \EE \FF_1 ) 
=  
\varepsilon^{-1} 
(\ii \omega(\vec{k}_0) + \ii \nabla\omega(\vec{k}_0  )\cdot (\vec{k}- \vec{k}_0) + \varepsilon^2 \partial_T) 
\hat{A}^{\varepsilon,per}_1(\frac{\vec{k}-\vec{k}_0}{\varepsilon},\varepsilon^2 t) 
\EE \FF_1 
$$
on the left-hand side of \eqref{syst4},
 expanding $ \omega(m \vec{k}_0 + \varepsilon \vec{K}) $ for $ m \in \{-3,-1,1,3\} $ 
w.r.t. $ \varepsilon $ on the right-hand side of \eqref{syst4}, 
in particular using 
$$
\omega(\vec{k}_0 + \varepsilon \vec{K}) = \omega(\vec{k}_0 ) + \varepsilon \nabla\omega(\vec{k}_0 ) \cdot  \vec{K} + \frac12 \varepsilon^2  \vec{K}^T\nabla^2 \omega(\vec{k}_0) \vec{K} + \mathcal{O}(\varepsilon^3),
$$
and equating the coefficients of 
{$ \varepsilon \EE$, $\varepsilon \EE^{-1}$, $ \varepsilon \EE^3$ and $ \varepsilon \EE^{-3}$} 
to zero 
yields with $ \vec{k} = m \vec{k}_0 + \varepsilon \vec{K} $  the NLS equation 
\begin{eqnarray*}
\partial_T \hat{A}_1(\vec{K} ,T) & = &  \frac12 \ii \vec{K}^T\nabla^2 \omega(\vec{k}_0) \vec{K} \hat{A}_1(\vec{K} ,T) + \frac{3   \omega_x^2(k_0) }{8 i \omega(\vec{k}_0)}
(\hat{A}_1 * \hat{A}_1*\hat{A}_{-1})(\vec{K} ,T) \\&&  - \frac{3(e^{ik_0}-1)(1-e^{-il_0})}{8 i \omega(\vec{k}_0)} (\hat{B}_1 * \hat{B}_1*\hat{B}_{-1})(\vec{K} ,T) 
\end{eqnarray*}
and the 
set of equations for the higher order corrections
\begin{eqnarray*}
- \ii \omega(\vec{k}_0) \hat{A}_{1,-1}(\vec{K} ,T)  & = & \ii \omega(-\vec{k}_0) \hat{A}_{1,-1}(\vec{K} ,T) 
+ \frac{3\omega_x^2(- k_0)}{8 i \omega(-\vec{k}_0)}
(\hat{A}_1 * \hat{A}_{-1}*\hat{A}_{-1})(\vec{K} ,T) \\ &&  - \frac{3(e^{-ik_0}-1)(1-e^{il_0})}{8 i \omega(-\vec{k}_0)} (\hat{B}_1 * \hat{B}_{-1}*\hat{B}_{-1})(\vec{K} ,T) ,
  \\
3 \ii \omega(\vec{k}_0) \hat{A}_{1,3}(\vec{K} ,T)  & = & \ii \omega(3 \vec{k}_0) \hat{A}_{1,3}(\vec{K} ,T) 
+ \frac{\omega_x^2(3 k_0)}{8 i \omega(3\vec{k}_0)}
(\hat{A}_1 * \hat{A}_1*\hat{A}_{1})(\vec{K} ,T) \\ &&  - \frac{(e^{3ik_0}-1)(1-e^{-3il_0})}{8 i \omega(3\vec{k}_0)} (\hat{B}_1 * \hat{B}_1*\hat{B}_{1})(\vec{K} ,T) ,
 \\
- 3 \ii \omega(\vec{k}_0) \hat{A}_{1,-3}(\vec{K} ,T)  & = & \ii \omega(-3 \vec{k}_0) \hat{A}_{1,-3}(\vec{K} ,T) 
+ \frac{\omega_x^2(-3k_0)}{8 i \omega(-3\vec{k}_0)}
(\hat{A}_{-1} * \hat{A}_{-1}*\hat{A}_{-1})(\vec{K} ,T) \\ &&  - \frac{(e^{-3ik_0}-1)(1-e^{3il_0})}{8 i \omega(-3\vec{k}_0)} (\hat{B}_{-1} * \hat{B}_{-1}*\hat{B}_{-1})(\vec{K} ,T) ,
\end{eqnarray*}
and associated equations for the $ \hat{A}_{-1} ,\ldots, \hat{A}_{-1,-3} $
and the $ \hat{B}_{1} ,\ldots, \hat{B}_{-1,-3} $.
Using the compatibility condition 
\eqref{UVrelation}, we obtain
at the wave vector $ \vec{k}_0 $ 
\begin{align} \label{comp_cond_Fourier_A_B} 
(e^{ik_0}-1)  \widehat{B}_1 (\vec{K} , T ) =  (e^{il_0}-1) \widehat{A}_1 ( \vec{K}, T) 
\end{align}
and, equivalently, due to our assumption $ A_{-1} = \overline {A_1} $ and $ B_{-1} = \overline {B_1} $, 
at the wave vector $- \vec{k}_0 $ 
$$
(e^{-ik_0}-1)  \widehat{B}_{-1} (\vec{K}, T)  =  (e^{-il_0}-1) \widehat{A}_{-1} ( \vec{K}, T) ,
$$ 
and the NLS equation in Fourier space obtained above becomes
\begin{align} \label{NLS_A_Fourier} 
\partial_T \hat{A}_1(K,T)  =   \frac12 \ii \vec{K}^T\nabla^2 \omega(\vec{k}_0) \vec{K} \hat{A}_1(K,T) + \gamma_A
(\hat{A}_1 * \hat{A}_1*\hat{A}_{-1})(K,T) ,
\end{align}
with 
\begin{equation} \label{gammaA}
\gamma_A =  \frac{3 \omega_x^2(k_0)}{8 i \omega(\vec{k}_0)} + \frac{3 \omega_y^4(l_0)}{8 i \omega_x^2(k_0) \omega(\vec{k}_0)} .
\end{equation}
It cannot be expected that the coefficient $ \gamma_A $ is symmetric 
in $ k_0 $ and $ l_0 $, since $ \hat{A}_1 $ describes the strain variable in $ x $-direction.
$ \hat{B}_1 $ describes the strain variable in $ y $-direction and, 
using \eqref{NLS_A_Fourier} and the compatibility condition \eqref{comp_cond_Fourier_A_B}, we obtain
that $ \hat{B}_1 $ satisfies 
\begin{align*} 
\partial_T \hat{B}_1(K,T)  =   \frac12 \ii \vec{K}^T\nabla^2 \omega(\vec{k}_0) \vec{K} \hat{B}_1(K,T) + \gamma_B
(\hat{B}_1 * \hat{B}_1*\hat{B}_{-1})(K,T) ,
\end{align*}
with 
$$ 
\gamma_B 
 = \gamma_A \frac{\omega_x^2(k_0)}{\omega_y^2(l_0)}
=  \frac{3 \omega_y^2(l_0)}{8 i \omega(\vec{k}_0)} + \frac{3 \omega_x^4(k_0)}{8 i \omega_y^2(l_0) \omega(\vec{k}_0)} .
$$ 
The  equation for $ \hat{B}_1 $ and the other $\hat{B}_{i,j}$ can also be  obtained directly as above, 
by inserting the corresponding ansatz $ \hat{\psi}_{V,1} $, $ \hat{\psi}_{V,-1} $ into \eqref{syst4}-\eqref{syst5} and using the compatibility condition  \eqref{comp_cond_Fourier_A_B}. 

Note that the chosen fixed wave vector $ \vec{k}_0 \in (\mathbb{T} \setminus \{ 0 \} )^2$, such that $ \omega(\vec{k}_0) > 0 $, allows us to determine $\hat A_{1,-1}$ by the formula given above. However, in order to determine $\hat A_{1,3}$, $\hat A_{1,-3}$ we need that moreover the non-resonance conditions 
\begin{equation} \label{nonres}
3 \omega(\vec{k}_0) \neq \omega(3\vec{k}_0)  \qquad \textrm{and} \qquad  \omega(3\vec{k}_0)  > 0
\end{equation}
hold true.
Note also that  \eqref{NLS_A_Fourier} yields \eqref{NLSintro}, 
with a factor $4$ in front of the coefficient $\gamma_A$, since the ansatz \eqref{uapprox1} corresponds to 
$\frac{1}{2} (\widehat{\psi}_{U,1} +\widehat{\psi}_{U,-1}) $ given by \eqref{eq21a}-\eqref{eq21b}, 
which yields 
$A = \frac{1}{2}A_1$.
\begin{remark} \label{remAB} 
In case $ k_0 \neq 0 $ and $ l_0 = 0 $ one must first derive the equation for $ \hat{A}_1$. 
Then, one obtains $ \hat{B}_1 \equiv 0 $ by  \eqref{comp_cond_Fourier_A_B} 
and subsequently $ \hat{B}_{-1}  \equiv 0 $ and $\hat{B}_{i,j} \equiv 0 $ 
by the analogous formulas to the ones given above for $ \hat{A}_{-1} $ and $\hat{A}_{i,j}$. 
This yields $ \hat{\psi}_{V,1}  \equiv \hat{\psi}_{V,-1}  \equiv 0 $ 
and corresponds to the 1D case of a modulated plane wave  
traveling in $x$-direction,  
cf.\ with (9) and (30) in \cite{Schn10}. 
Analogously this holds true in case $ l_0 \neq 0 $ and $ k_0 = 0 $. 
\end{remark}

v) By this construction, 
 that is, 
by inserting $\hat{\psi}_{U,\pm 1}(\vec{k},t)$, $\hat{\psi}_{V,\pm1}(\vec{k},t)$   into \eqref{syst4}-\eqref{syst5} 
with  $\hat A$, $\hat B$ determined by the NLS equation \eqref{NLS_A_Fourier} for $\hat A_1$ 
and the compatibility condition \eqref{comp_cond_Fourier_A_B},
we eliminated all terms of  order 
$ \mathcal{O}(\varepsilon^3) $, if we measure the magnitude in $ L^1 $. 
Using 
\begin{equation*}
\omega(\vec{k}) = \omega(\vec{k}_0 ) + \nabla\omega(\vec{k}_0 ) \cdot  (\vec{k}-\vec{k}_0) 
+ \frac12 (\vec{k}-\vec{k}_0)^T\nabla^2 \omega(\vec{k}_0) (\vec{k}-\vec{k}_0) 
+ R_{\omega, \vec{k}_0} (\vec{k}-\vec{k}_0) ,
\end{equation*} 
we obtain 
\begin{eqnarray*}
&& \|  R_{\omega, \vec{k}_0} (\vec{k}-\vec{k}_0) 
\varepsilon^{-1} \hat{A}^{\varepsilon,per}_1(\frac{\vec{k}-\vec{k}_0}{\varepsilon},\varepsilon^2 t) \|_{L^1(\mathbb{T}^2,d\vec{k})}
\leq  C \varepsilon^4 \|\hat{A}^{\varepsilon,per}_1(\cdot,T) \|_{L^1_3(\mathbb{T}_{2 \pi/\varepsilon}^2)}
\\ && 
\leq C \varepsilon^4 \|\hat{A}_1(\cdot,T) \|_{L^1_3(\R^2)} 
\leq C \varepsilon^4 \|\hat{A}_1(\cdot,T) \|_{L^2_s(\R^2)} 
\leq C \varepsilon^4 \| A_1(\cdot,T) \|_{H^s(\R^2)}, 
\quad s > 4, 
\end{eqnarray*}
where we used the Sobolev embedding $ L^2_s(\R^2) \subset L^1_3(\R^2) $
for $ s > 4 $ and the fact that the Fourier transform is an isomorphism between 
$ L^2_s(\R^2) $ and $ H^s(\R^2) $. 
Recalling that the $L^1(\mathbb{T}^2)$-norm is closed under convolution and noting that the other terms in \eqref{eq21a}-\eqref{eq21b} have a higher order by a factor $ \mathcal{O}(\varepsilon^2) $,
these and similar estimates 
show that the remaining terms which are collected in the residual are of order $ \mathcal{O}(\varepsilon^4) $ in $ L^1 $. 
Moreover, we have 
$$ 
\| \widehat{A} -  \widehat{A}^{\varepsilon} \|_{L^1(\R^2)} \leq C  \varepsilon^3 \| \widehat{A} \|_{L^1_3 (\R^2)} 
$$ 
such that $ \widehat{\psi}_{u,m,n}  - \frac{1}{2} (\widehat{\psi}_{U,1} +\widehat{\psi}_{U,-1} ) $ are  $ \mathcal{O}(\varepsilon^3) $-close  in $ L^1 $ and, thus, their counterparts in physical space are $ \mathcal{O}(\varepsilon^3) $-close  in $ L^{\infty} $,  due to 
the well known inequality
\begin{equation} \label{l1linf}
\sup_{m,n \in \Z^2}|r_{m,n}| \leq  
\sup_{m,n \in \Z^2} \int_{\mathbb{T}^2} |\hat{r}(k,l) e^{ikm+iln}| dk dl
\leq \int_{\mathbb{T}^2} |\hat{r}(k,l) | dk dl = \| \hat{r} \|_{L^1}.
\end{equation}
Also note  that $ \| \widehat{\psi}_{U,\pm 1} \|_{L^1(\mathbb{T}^2)} = \mathcal{O}(\varepsilon) $. This will be used in the next section.\\
vi) Finally, 
our approximation ansatz 
has to be slightly modified to satisfy the compatibility
condition \eqref{UVrelation}. We define the space of functions satisfying the compatibility equations
$$ 
\mathcal{X} = \{ U : = (\widehat{U}_1,\widehat{U}_{-1},\widehat{V}_1,\widehat{V}_{-1}): \mathbb{T}^2 \to \C^4 :  
(e^{ik}-1)  \widehat{V}_{\pm1}(k,l,t)  =  (e^{il}-1) \widehat{U}_{\pm1}(k,l,t) \},
 $$ 
equipped with the norm of $(L^1)^4$,
$$ 
\|U \|_{\mathcal{X}} = \|\widehat{U}_1 \|_{L^1} + \|\widehat{U}_{-1} \|_{L^1}+ \|\widehat{V}_1 \|_{L^1}+ \|\widehat{V}_{-1}\|_{L^1}.
$$ 
%
 Setting for shortness $a= e^{ik}-1$, $b= e^{il}-1$, 
the projection $ \widehat{P} $ from $ (L^1)^4 $ onto $ \mathcal{X} $ is given by
$$
\widehat{P} U  
= 
\frac{1}{a^2+b^2}
( 
a ( a \widehat{U}_1 + b \widehat{V}_1 ) , 
a ( a \widehat{U}_{-1} + b \widehat{V}_{-1} ) ,
b ( a \widehat{U}_1 + b \widehat{V}_1 ) , 
b ( a \widehat{U}_{-1} + b \widehat{V}_{-1} ) 
)^T  , 
$$ 
such that 
$$
U - \widehat{P} U  
= 
\frac{1}{a^2+b^2}
( b (b \widehat{U}_1 - a \widehat{V}_1 ) , 
b ( b \widehat{U}_{-1}  - a \widehat{V}_{-1} ) ,
a ( a \widehat{V}_1  -  b \widehat{U}_1 )  , 
a ( a \widehat{V}_{-1} - b \widehat{U}_{-1} ) 
)^T  . 
$$ 
We set $ \varepsilon \Psi_{U,V} =  (\hat{\psi}_{U,1} , \hat{\psi}_{U,-1} , \hat{\psi}_{V,1} , \hat{\psi}_{V,-1}) $
and $ \varepsilon \Psi = \widehat{P} \varepsilon \Psi_{U,V}  $.
Thus, 
recalling \eqref{comp_cond_Fourier_A_B},  the definition of $\hat{\psi}_{U,\pm1}$, $\hat{\psi}_{V,\pm1}$,
and using an analogous argument as for the first estimate in v) above, we obtain
$$ 
\| \varepsilon \Psi -  \varepsilon \Psi_{U,V} \|_{(L^1)^4} \leq C \varepsilon^2.
$$ 
Moreover, since for initial data of the original system which satisfy the compatibility condition its solutions stay in $\mathcal{X}$ for all times, we obtain  by the linearity and boundedness of the projection operator $\widehat{P} $ in 
$(L^1)^4$ 
that  
by 
inserting $ \varepsilon  \Psi  $ into \eqref{syst4}-\eqref{syst5} 
the residual
is of order $ \mathcal{O}(\varepsilon^4) $ in $(L^1)^4$  and belongs to $ \mathcal{X} $,
i.e., also satisfies the compatibility condition \eqref{UVrelation}.

\section{The error estimates}
\label{sec4}

The system \eqref{syst4}-\eqref{syst5} in Fourier space is abbreviated by 
$$ 
\partial_t U = L U + C(U,U,U),
$$
where $ L $ stands for the linear terms and where $ C $ is a symmetric trilinear mapping.
The operator $ L $ defines a uniformly bounded semigroup $ (e^{Lt})_{t \geq 0} $ in the phase space $ \mathcal{X} \subset (L^1(\mathbb{T}^2))^4 $, i.e., there exists a $ C_L > 0 $, here
$ C_L = 1 $, such that $ \sup_{t \geq 0} \| e^{Lt} \|_{\mathcal{X} \to \mathcal{X}}  \leq C_L$.
The trilinear mapping $ C $ satisfies $ \| C(U,V,W) \|_{\mathcal{X}} \leq C_C \| U \|_{\mathcal{X}}  \| V \|_{\mathcal{X}}
 \| W \|_{\mathcal{X}} $ with a constant $ C_C $ independent of $ U,V,W \in \mathcal{X} $.

The error $ \varepsilon^2 R = U - \varepsilon \Psi $ made by the approximation $  \varepsilon \Psi $ from Section \ref{sec3} 
 satisfies 
$$ 
\partial_t R = L R + 3 \varepsilon^2 C(\Psi,\Psi,R) + 
3 \varepsilon^{3} C(\Psi,R,R) + \varepsilon^{4} C(R,R,R)
 +  \varepsilon^{-2} \textrm{Res}(\varepsilon \Psi),
$$ 
where the residual 
$$
\textrm{Res}(\varepsilon \Psi) = - \partial_t \varepsilon \Psi + L \varepsilon \Psi + C(\varepsilon \Psi,\varepsilon \Psi,\varepsilon \Psi),
$$
contains all terms which do not cancel after inserting the approximation $ \varepsilon \Psi $ 
into the equations. 
By the above construction we have 
\begin{lemma}
For the approximation $ \varepsilon \Psi $  there exist $ \varepsilon_0 > 0 $  and  $ C_{res} > 0 $ such that  
for $ \varepsilon \in (0,\varepsilon_0) $ we have
$$
\sup_{t \in [0,T_0/\varepsilon^2]}  \|\varepsilon^{-2} \textrm{Res}(\varepsilon \Psi)\|_{\mathcal{X}}  \leq  C_{res}  \varepsilon^2 .
$$
\end{lemma}
In order to estimate the error function $ R $ on the time interval $  [0,T_0/\varepsilon^2] $ we rewrite the 
equation for $ R $ with the help of the variation of constant formula as 
$$ 
R(t) = e^{Lt} R(0)  + \int_0^t e^{L(t-\tau)}(3 \varepsilon^2 C(\Psi,\Psi,R) + 
3 \varepsilon^{3} C(\Psi,R,R) + \varepsilon^{4} C(R,R,R)
+ \varepsilon^{-2} \textrm{Res}(\varepsilon \Psi))(\tau) d \tau.
$$ 
We estimate 
$$ 
\|R(t) \| _{\mathcal{X}} \leq C_L  \|R(0) \| _{\mathcal{X}} + \int_0^t C_L (3  \varepsilon^2 C_C C_{\Psi}^2 \|R(\tau) \| _{\mathcal{X}}  + 3  \varepsilon^3 C_C C_{\Psi} \|R(\tau) \| _{\mathcal{X}}^2 
+  \varepsilon^4 C_C \|R(\tau) \| _{\mathcal{X}}^3
+ C_{res}  \varepsilon^2 )  d\tau
$$ 
where $ C_{\Psi} = \sup_{t \in  [0,T_0/\varepsilon^2]} \|\Psi(t) \| _{\mathcal{X}}$.
As long as
$$ 3  \varepsilon C_C C_{\Psi} \|R(\tau) \| _{\mathcal{X}} 
+  \varepsilon^2 C_C \|R(\tau) \| _{\mathcal{X}}^2 \leq 1  $$
holds true, we have
$$ 
\|R(t) \| _{\mathcal{X}} \leq C_L  \|R(0) \| _{\mathcal{X}} + \int_0^t C_L (  \varepsilon^2 (3 C_C C_{\Psi}^2 + 1)\|R(\tau) \| _{\mathcal{X}} 
+ C_{res}  \varepsilon^2 )  d\tau.
$$ 
Gronwall's inequality gives 
$$
\|R(t) \| _{\mathcal{X}} \leq C_L  (\|R(0) \| _{\mathcal{X}} +  C_{res} T_0) e^{C_L   (3 C_C C_{\Psi}^2 + 1) T_0} =: M
$$
for all $ t \in [0,T_0/\varepsilon^2] $. Choosing $ \varepsilon_0 > 0 $ so small that 
$ 
3  \varepsilon_0 C_C C_{\Psi} M
+  \varepsilon_0^2 C_C M^2 \leq 1  
$
is satisfied, we are done in Fourier space. The estimate stated in 
Theorem \ref{th11} follows by undoing the above transformations and applying the  inequality \eqref{l1linf}
to the error $ R (t) \in  \mathcal{X} = (L^1(\mathbb{T}^2))^4 $ in Fourier space.


\section{Discussion}
\label{secdisc}

It is the purpose of this last section to discuss two additional topics, namely how to obtain 
an approximation result for the original displacement variables $ q_{m,n} $  and secondly the 
robustness of 
this
result 
w.r.t. small variations of the interaction force $ W' $. 
%

\subsection{The approximation result for the displacement variables $ q_{m,n} $}

In the following we would like to explain that an NLS approximation result for the displacement variables $ q_{m,n} $
can be obtained without  requiring
an on-site potential as used so far in the existing literature, 
cf. \cite{GM04} for the one-dimensional and \cite[Theorem 7.1 with (7.1)]{GHM06} for the multi-dimensional case,  
where the results are derived in physical space and justified with respect to the $\ell^2$-energy norm.  
%
Note that due to scaling, the estimates in the energy norm imply that in the multidimensional case one has to include also higher-order corrections in the approximation. In the presence of an on-site potential we have $ \omega(0) \neq 0 $ which simplifies the derivation and justification of the NLS equation in contrast to the case without on-site potential where we have $ \omega(0) = 0 $. The presence of a non-oscillating term corresponding to $ \omega(0)=0 $ may yield in principle a more complex structure of modulation equations. 

Concerning the approach taken here, we remark moreover that 
from 
Theorem \ref{th11}  we would obtain a different non-satisfying approximation result due to the required reconstruction of the displacement 
variables from the strain variables. Therefore, we proceed as follows.

The FPUT system \eqref{fpuintro}  with the interaction potential \eqref{Wprime} is given by 
\begin{eqnarray}   \label{qmneq}
\partial_t^2 q_{m,n} & = & q_{m+1,n}- 2 q_{m,n} -q_{m-1,n} - (q_{m+1,n}-q_{m,n})^3 + (q_{m,n}-q_{m-1,n})^3 
\\ && 
+ q_{m,n+1}-2 q_{m,n} -q_{m,n-1}
- (q_{m,n+1}-q_{m,n})^3 + (q_{m,n}-q_{m,n-1})^3 . \nonumber
\nonumber
\end{eqnarray}
With the same $\omega$ as in \eqref{linear-pr} and the 
periodic convolution \eqref{perconv}
this system reads in Fourier space as 
\begin{align}  
\partial_t^2 \hat{q}(k,l,t)   
=  - \omega^2(k,l) \hat{q}(k,l,t) 
- \hat{N} ( \hat{q})(k,l,t) , 
\label{fpufourier}
\end{align}
where 
\begin{eqnarray*}
\hat{N} ( \hat{q})(k,l,t) 
& = & ( ( e^{ik} - 1 ) \hat{q}(k,l,t) )^{\ast 3} - ( ( 1 - e^{-ik} ) \hat{q}(k,l,t)  )^{\ast 3} 
\\&& \qquad \qquad
+  ( ( e^{il} - 1 ) \hat{q}(k,l,t) )^{\ast 3} - ( ( 1 - e^{-il} ) \hat{q}(k,l,t)  )^{\ast 3}
\\ 
& = &
\int_{\mathbb{T}^2} \int_{\mathbb{T}^2}
( n (k-k_1, k_1-k_2,k_2 ) + n ( l - l_1  , l_1- l_2, l_2) ) 
\\ && \qquad \qquad \times 
\hat{q}(k-k_1, l - l_1  ,t) 
\hat{q}(k_1-k_2,l_1- l_2,t) \hat{q}(k_2,l_2,t) d (k_2,l_2) d (k_1,l_1) ,
\end{eqnarray*}
with 
$$ 
n (k_1,k_2,k_3) = ( e^{ik_1} - 1 ) ( e^{ik_2} - 1 ) ( e^{ik_3} - 1 )  + c.c. .
$$
Note that
$$ 
n (k_1,k_2,k_3) =n (k_1,k_3,k_2) =n (k_2,k_1,k_3) =  n (k_2,k_3,k_1) = n (k_3,k_1,k_2)  = n (k_3,k_2,k_1) . 
$$
Using the abbreviations $ \vec k = (k,l) $, $ \vec k_i = (k_i,l_i) $ for  $ i = 1,2, 3 $, 
and 
\begin{equation} \label{D}
D ( \vec k_1, \vec k_2 , \vec k_3 ) = n ( k_1,k_2 , k_3 ) + n ( l_1, l_2, l_3) 
\end{equation} 
we obtain the more concise formulation
$$
\hat{N} ( \hat{q})(\vec k,t) 
 = 
\int_{\mathbb{T}^2} \int_{\mathbb{T}^2}
D ( \vec k  - \vec k_1  , \vec k_1 - \vec k_2 , \vec k_2) 
\hat{q}(\vec k - \vec k_1   ,t) 
\hat{q}(\vec k_1- \vec k_2 ,t) \hat{q}(\vec k_2,t) d \vec k_2  d \vec k_1.
$$
The above kernel $ n $ can be estimated  as follows 
\begin{lemma} 
There exists a $ C > 0 $ such that for 
$ k_1, k_2, k_3 \in \mathbb T $
we have 
$$
| n (k_1,k_2,k_3)   
|_{k =k_1 +k_2 + k_3} | \leq C |k| .
$$
\end{lemma}
\noindent
{\bf Proof.} 
An elementary calculation shows 
\begin{eqnarray*}
n (k_1,k_2,k_3)  |_{k =k_1 +k_2 + k_3} 
& = &2 (\cos k - 1 ) ( 1 - \cos k_1   -  \cos k_2  -  \cos k_3 )|_{k =k_1 +k_2 + k_3} 
\\ && \qquad  \qquad 
-  2 (\sin k)( \sin k_1  + \sin k_2  + \sin k_3  )|_{k =k_1 +k_2 + k_3} ,
\end{eqnarray*}
which can be estimated by 
$$
| n (k_1,k_2,k_3) |_{k =k_1 +k_2 + k_3} |
\leq 8 |\cos k -1|  + 6 |\sin k| \leq C |k| ,
$$
with a constant $ C $ independent  of $ k $.
\qed 
\medskip
%

\begin{remark}{\rm
Note that the limit $  \lim_{k\to0} n (k_1,k_2,k_3) |_{k =k_1 +k_2 + k_3}/k $ does not exist.  
Compare, e.g., $k_1=k_2=k_3 = 1/n$ with $k_1=k_2=\pi/4 $, $k_3=-(\pi/2) +(1/n) $ for $n\to \infty$.      
}\end{remark}

This property allows us to rewrite \eqref{fpufourier} as a first order system.
Indeed, since $ \omega (k,l) = \mathcal {O} ( (k^2 + l^2)^{1/2} ) $ we have 
$$
\sup_{ { k,l \in \mathbb T} } | (i \omega(k,l))^{-1}  k | \leq C < \infty , 
\qquad  \textrm{and} \qquad 
\sup_{{ k,l \in \mathbb T} } | (i \omega(k,l))^{-1}  l | \leq C < \infty , 
$$
and 
we find 
\begin{eqnarray}  
\label{fpufourier1}
\partial_t \hat{q} & = & i \omega  \hat{q}_1 , \\
\partial_t \hat{q}_1& = & i \omega  \hat{q}  
- 
\frac1{i \omega}
\hat{N} (\hat{q})  ,
\label{fpufourier2}
\end{eqnarray}
where $ \frac1{i \omega}
\hat{N} (\hat{q})  $ is a smooth mapping from $ L^1 $ to $ L^1 $. Therefore, we have a system  
with the same properties as in the previous sections and so the same 
analysis as above can be carried out  for  \eqref{fpufourier1}-\eqref{fpufourier2}.
The system is diagonalized  by introducing $\hat Q_1 = \hat q + \hat q_1$ and  $\hat Q_{-1} = \hat q - \hat q_1$
such that 
\begin{eqnarray}  
\label{fpufourier1diagonal}
\partial_t \hat{Q}_1 & = & i \omega  \hat{Q}_1 
- \frac1{8 i \omega} 
 \hat{N}  (\hat Q_1 + \hat Q_{-1}) 
 , \\
\partial_t \hat{Q}_{-1}& = & - i \omega  \hat{Q}_{-1}   
+ 
\frac1{8 i \omega}
 \hat{N}  (\hat Q_1 + \hat Q_{-1}) 
\label{fpufourier2diagonal}
\end{eqnarray}
For $\hat Q_1$ and  $\hat Q_{-1}$ we make   
the same ansatz \eqref{eq21a}-\eqref{eq21b} as for 
$\hat U_1$ and  $\hat U_{-1}$
in \eqref{syst4}-\eqref{syst5} 
except that now the ansatz is  called $\hat{\psi}_{Q,1}$ and  $\hat{\psi}_{Q,-1} $ in the 
following. 
The derivation of the equations to obtain $A_1 , \ldots, A_{-1,-3}$ follows the same pattern as in the case of the strains
discussed in Section \ref{sec3}. 
The only significant difference arises  from the nonlinear terms 
$ \hat{N}( \hat Q_1+ \hat Q_{-1} ) $. 
At order $\varepsilon$ we get 
\begin{align*}
& \varepsilon^{-3} \EE^3 \FF_3 
\int_{\mathbb{T}^2} \int_{\mathbb{T}^2}
D ( \vec k  - \vec k_1  , \vec k_1 - \vec k_2, \vec k_2) 
\hat{A}^{\varepsilon,per}_1(\frac{ \vec k - \vec k_1 -\vec{k}_0}{\varepsilon},T)
\\ & \qquad \qquad \qquad \qquad  \times 
\hat{A}^{\varepsilon,per}_1(\frac{ \vec k_1- \vec k_2   -\vec{k}_0}{\varepsilon},T)
\hat{A}^{\varepsilon,per}_1(\frac{ \vec k_2 - \vec{k}_0}{\varepsilon},T)
d \vec k_1  d \vec k_2
\\
& +3 \varepsilon^{-3}  \EE \FF_1 
\int_{\mathbb{T}^2} \int_{\mathbb{T}^2}
D ( \vec k  - \vec k_1  , \vec k_1 - \vec k_2, \vec k_2) 
\hat{A}^{\varepsilon,per}_1(\frac{ \vec k - \vec k_1  -\vec{k}_0}{\varepsilon}, T)
\\ & \qquad \qquad \qquad \qquad  \times 
\hat{A}^{\varepsilon,per}_1(\frac{ \vec k_1 - \vec k_2-\vec{k}_0}{\varepsilon},T)
\hat{A}^{\varepsilon,per}_{-1}(\frac{ \vec k_2 + \vec{k}_0}{\varepsilon},T)
d \vec k_1  d \vec k_2
\\
& + 3 \varepsilon^{-3}  \EE^{-1} \FF_{-1}
\int_{\mathbb{T}^2} \int_{\mathbb{T}^2}
D ( \vec k  - \vec k_1  , \vec k_1- \vec k_2, \vec k_2) 
\hat{A}^{\varepsilon,per}_1 (\frac{ \vec k - \vec k_1 -\vec{k}_0}{\varepsilon},T)
\\ & \qquad \qquad \qquad \qquad  \times 
\hat{A}^{\varepsilon,per}_{-1}(\frac{ \vec k_1 - \vec k_2+ \vec{k}_0}{\varepsilon},T)
\hat{A}^{\varepsilon,per}_{-1}(\frac{ \vec k_2 + \vec{k}_0}{\varepsilon},T)
d \vec k_1  d \vec k_2
\\
& + \varepsilon^{-3} \EE^{-3} \FF_{-3} 
\int_{\mathbb{T}^2} \int_{\mathbb{T}^2}
D ( \vec k  - \vec k_1  , \vec k_1- \vec k_2, \vec k_2) 
\hat{A}^{\varepsilon,per}_{-1}(\frac{ \vec k - \vec k_1 + \vec{k}_0}{\varepsilon},T)
\\ & \qquad \qquad \qquad  \qquad \times 
\hat{A}^{\varepsilon,per}_{-1}(\frac{ \vec k_1 - \vec k_2 + \vec{k}_0}{\varepsilon},T)
\hat{A}^{\varepsilon,per}_{-1}(\frac{ \vec k_2 + \vec{k}_0}{\varepsilon},T)
d \vec k_1  d \vec k_2.
\end{align*}
As in Section \ref{sec3}, 
we set, correspondingly, 
$ \varepsilon \vec K_1 = \vec k_1 - \vec k_2 \pm \vec k_0 $, $ \varepsilon  \vec K_2 = \vec k_2 \pm \vec k_0 $ 
and $\varepsilon \vec K = \vec k - m\vec k_0 $, $m \in \{ - 3, -1, 1, 3 \}$ 
and take the limit  $ \varepsilon \to 0 $ to obtain 
\begin{align*}
%
&{  \varepsilon  \EE^3 }
\int_{\R^2} \int_{\R^2}
D (\vec k_0  , \vec{k}_0  , \vec{k}_0  ) 
\hat{A}_1( \vec K - \vec K_1 - \vec K_2 ,T)
\hat{A}_1( \vec K_1 ,T)
\hat{A}_1( \vec K_2 ,T)
d \vec K_1  d \vec K_2
\\
& + 3 {  \varepsilon   \EE }
\int_{\R^2} \int_{\R^2}
D ( \vec k_0 , \vec k_0  , - \vec k_0 ) 
\hat{A}_1( \vec K  - \vec K_1 - \vec K_2 ,T)
\hat{A}_1( \vec K_1 ,T)
\hat{A}_{-1}( \vec K_2 ,T)
d \vec K_1  d \vec K_2
\\
& + 3{ \varepsilon   \EE^{-1} } 
\int_{\R^2} \int_{\R^2}
D ( \vec k_0  , -  \vec k_0 , - \vec k_0 ) 
\hat{A}_1 ( \vec K - \vec K_1 - \vec K_2 ,T)
\hat{A}_{-1}( \vec K_1 ,T)
\hat{A}_{-1}( \vec K_2 ,T)
d \vec K_1  d \vec K_2
\\
& + { \varepsilon \EE^{-3} } 
\int_{\R^2} \int_{\R^2}
D ( - \vec k_0  , - \vec k_0  , - \vec k_0  ) 
\hat{A}_{-1}( \vec K  - \vec K_1 - \vec K_2 ,T)
\hat{A}_{-1}( \vec K_1 ,T)
\hat{A}_{-1}( \vec K_2 ,T)
d \vec K_1  d \vec K_2.
\end{align*}

Thus, as in Section \ref{sec3}, we  obtain  the NLS equation in Fourier space
\begin{equation} \label{NLS5}
\partial_T \hat{A}_1(K,T)  =   \frac12 \ii \vec{K}^T\nabla^2 \omega(\vec{k}_0) \vec{K} \hat{A}_1(K,T) 
- \frac{3 D (\vec k_0  , \vec{k}_0  , -\vec{k}_0 ) }{8 i \omega(\vec{k}_0)} (\hat{A}_1 * \hat{A}_1*\hat{A}_{-1})(K,T) ,
\end{equation}
at
${ \varepsilon  \EE } $,
where, by \eqref{omega} and \eqref{D}, 
$ D ( \vec k_0, \vec k_0 , -\vec k_0 ) = - (\omega_x^4(k_0) + \omega_y^4(l_0) )$, 
and the equations 
\begin{eqnarray*}
- \ii \omega(\vec{k}_0) \hat{A}_{1,-1}(K,T)  & = & \ii \omega(-\vec{k}_0) \hat{A}_{1,-1}(K,T) 
- \frac{3 D (\vec k_0  , - \vec{k}_0  , - \vec{k}_0  ) }{8 i \omega(-\vec{k}_0)}(\hat{A}_1 * \hat{A}_{-1}*\hat{A}_{-1})(K,T) ,
  \\
3 \ii \omega(\vec{k}_0) \hat{A}_{1,3}(K,T)  & = & \ii \omega(3 \vec{k}_0) \hat{A}_{1,3}(K,T) 
- \frac{D (\vec k_0  , \vec{k}_0  , \vec{k}_0  )}{8 i \omega(3\vec{k}_0)} (\hat{A}_1 * \hat{A}_1*\hat{A}_{1})(K,T) ,
 \\
- 3 \ii \omega(\vec{k}_0) \hat{A}_{1,-3}(K,T)  & = & \ii \omega(-3 \vec{k}_0) \hat{A}_{1,-3}(K,T) 
- \frac{D ( - \vec k_0  , - \vec{k}_0  , - \vec{k}_0  ) }{8 i \omega(-3\vec{k}_0)} (\hat{A}_{-1} * \hat{A}_{-1}*\hat{A}_{-1})(K,T) ,
\end{eqnarray*}
for the higher order corrections at  
$ \varepsilon \EE^{-1} $, $ \varepsilon  \EE^3 $, $ \varepsilon \EE^{-3} $,
and associated equations for the $ \hat{A}_{-1} ,\ldots, \hat{A}_{-1,-3} $. 
As in Section \ref{sec3}, by this construction the residual terms are  of order 
$ \mathcal{O}(\varepsilon^4) $ in $ (L^1)^2 $. (Of course, when considering directly the displacements, there is no compatibility condition to be satisfied.) 

Thus, we are exactly in 
the same situation as at the beginning of Section \ref{sec4}.
The system \eqref{fpufourier1diagonal}-\eqref{fpufourier2diagonal} is abbreviated as 
\begin{equation} \label{Qeq} 
\partial_t Q = L Q + C(Q,Q,Q),
\end{equation}
where $ L $ stands for the linear terms and where $ C $ is a symmetric trilinear mapping.
The operator $ L $ defines a uniformly bounded semigroup $ (e^{Lt})_{t \geq 0} $ in the phase space $ \mathcal{X} = (L^1)^2 $, i.e., there exists a $ C_L > 0 $, here
$ C_L = 1 $, such that $ \sup_{t \geq 0} \| e^{Lt} \|_{\mathcal{X} \to \mathcal{X}}  \leq C_L$.
The trilinear mapping $ C $ satisfies $ \| C(U,V,W) \|_{\mathcal{X}} \leq C_C \| U \|_{\mathcal{X}}  \| V \|_{\mathcal{X}}
 \| W \|_{\mathcal{X}} $ with a constant $ C_C $ independent of $ U,V,W \in \mathcal{X} $.
The error $ \varepsilon^2 R = Q - \varepsilon \Psi $ made by the approximation $  \varepsilon \Psi = (\hat{\psi}_{Q,1},\hat{\psi}_{Q,-1}) $ 
 satisfies 
$$ 
\partial_t R = L R + 3 \varepsilon^2 C(\Psi,\Psi,R) + 
3 \varepsilon^{3} C(\Psi,R,R) + \varepsilon^{4} C(R,R,R)
+ \varepsilon^{-2} \textrm{Res}(\varepsilon \Psi).
$$ 
Following the rest of Section \ref{sec4} line for line  shows the $ \mathcal{O}(1) $-boundedness of $ \sup_{t \in [0,T_0/\varepsilon^2]} \| R(t) \|_{\mathcal{X}} $.
Thus, we have proved 
\begin{theorem} \label{th51}
Let $ s_A > 4 $ 
and $ (k_0,l_0) \neq (0,0) $ chosen in such a way that the non-resonance condition \eqref{nonres} is satisfied.
For all $ T_0 > 0 $, $ C_1 > 0 $, $ C_2 > 0 $  
there exist $ \varepsilon_0 > 0 $ and $ C_3 > 0 $ such that for all 
$ \varepsilon \in (0,\varepsilon_0) $ the following holds.
Let $ A \in C([0,T_0],H^{s_A}(\R^2,\C)) $ be a solution of the  NLS equation 
$$
\partial_T A  =  -  \frac12 \ii (\partial_X,\partial_Y) \nabla^2 \omega(\vec{k}_0) (\partial_X,\partial_Y)^T A  
- \frac{3 i (\omega_x^4(k_0) +\omega_y^4(l_0) )}{ 2 \omega(\vec{k}_0)} 
|A|^2 A,
$$
with 
$$ 
\sup_{T \in [0,T_0]} \| A(\cdot,T) \|_{H^{s_A}} \leq C_1 
$$ 
and let 
$$
\psi_{q,m,n} =  \varepsilon A(X,Y,T) e^{i (k_0 m + l_0 n + \omega_0 t) } + c.c.,
$$
with $X,Y,T$ as in \eqref{XYT}.
Take initial conditions of \eqref{qmneq}
with 
$$
\sup_{(m,n) \in \Z^2} (| q_{m,n}(0) -  \psi_{q,m,n}
(0) | +| \partial_t q_{m,n}(0) -  \partial_t \psi_{q,m,n}
(0) | ) \leq  C_2 \varepsilon^2.
$$
Then the  
solutions $ (q_{m,n})_{(m,n) \in \Z^2} $ of 
\eqref{qmneq} with these initial conditions  satisfy
$$
 \sup_{t \in [0,T_0/\varepsilon^2]} \sup_{(m,n) \in \Z^2} (| q_{m,n}(t) -  \psi_{q,m,n}
(t) | +| \partial_t q_{m,n}(t) -  \partial_t \psi_{q,m,n}
(t) |
 ) \leq C_3 \varepsilon^2.
$$ 
\end{theorem}
\begin{remark} \label{remABdisplacements} 
In analogy to Remark \ref{remAB}, the case $ k_0 \neq 0 $ and $ l_0 = 0 $
corresponds to the 1D case of a modulated plane wave traveling in $x$-direction,
cf.\ with (2.14) in \cite{GM04}.
Moreover, similar results are expected to hold true also for multi-dimensional lattices and general cubic interaction potentials with neighbors at an arbitrary (finite) distance, like the ones used for instance in \cite{G10}.
However, here we do not pursue this further. 
\end{remark}

\subsection{Small variations of the interaction force $ W' $}

We would like to discuss the robustness of our result w.r.t.\ small perturbations of the interaction forces.  
Therefore, we 
consider the FPUT system
\begin{eqnarray}  
\label{fpusec52}
\partial_t^2 q_{m,n} & = & W_{m \to m+1,n}'(q_{m+1,n}-q_{m,n})- W_{m-1 \to m,n}'(q_{m,n}-q_{m-1,n})\\ && +W_{m ,n\to n+1}'(q_{m,n+1}-q_{m,n})- W_{m ,n-1\to n}'(q_{m,n}-q_{m,n-1})
\nonumber
\end{eqnarray}
for all $ (m,n) \in \Z^2 $ where the interaction forces 
are small perturbations of the original 
interaction force, i.e., we consider 
$$ 
W_{m-1 \to m,n}'(u) = u + \alpha_{m-1 \to m,n} \varepsilon^3 u
+ \beta_{m-1 \to m,n} \varepsilon^2 u^2
-u^3 +  \gamma_{m-1 \to m,n} \varepsilon u^3
+ \mathcal{O}(u^4)
$$
and 
$$ 
W_{m,n \to n+1}'(u) = u + \alpha_{m,n \to n+1} \varepsilon^3 u
+ \beta_{m,n \to n+1} \varepsilon^2 u^2
-u^3 +  \gamma_{m,n \to n+1} \varepsilon u^3
+ \mathcal{O}(u^4).
$$
Assuming that the Fourier transforms of  $ (\alpha_{m-1 \to m,n})_{m,n \in \Z^2} ,\ldots , (\gamma_{m,n \to n+1})_{m,n \in \Z^2}$ are $ \mathcal{O}(1) $ in $ L^1 $, these additional terms do not affect the derivation of the NLS equation and
the equations for the higher order corrections, since after inserting $ \varepsilon \Psi = (\hat{\psi}_{Q,1},\hat{\psi}_{Q,-1}) 
$ into \eqref{fpufourier1diagonal}-\eqref{fpufourier2diagonal} the created new terms are of order $ \mathcal{O}(\varepsilon^4) $ in $ (L^1)^2 $ and add to the residual.
The counterpart to \eqref{Qeq} is of the form 
\begin{equation} \label{Qeq1} 
\partial_t Q = L Q + \varepsilon^3 P_1(Q) + \varepsilon^2 P_2(Q,Q) + C(Q,Q,Q) + \varepsilon P_3 (Q,Q,Q) + P_4(Q),
\end{equation}
with $ P_1 $ linear, $ P_2 $ bilinear, $ P_3 $ trilinear in their arguments and $ \| P_4(Q) \|_{\mathcal{X}}  
\leq C \| Q \|_{\mathcal{X}} ^4 $. Thus inserting 
$ Q = \varepsilon \Psi + \varepsilon^2 R $ in \eqref{Qeq1} yields the equation for the error 
\begin{equation} \label{Req1} 
\partial_t R = L R + 3 \varepsilon^2 C(\Psi,\Psi,R) + 
3 \varepsilon^{3} C(\Psi,R,R) + \varepsilon^{4} C(R,R,R)
+ \varepsilon^{-2} \textrm{Res}(\varepsilon \Psi) + \varepsilon^3 G(R) , 
\end{equation}
with $ \varepsilon^3 G(\cdot) $ a smooth mapping  in $ \mathcal{X}$ coming from the new terms $ P_1, \ldots, P_4 $.
Since $ \varepsilon^3 G(R)  $ is of order $ \mathcal{O}(\varepsilon^3) $, a straightforward modification of 
the proof given in Section \ref{sec4} allows to 
conclude
\begin{corollary}
Under  the above assumptions on the interaction forces  
Theorem \ref{th51} remains valid.
\end{corollary}

\bibliographystyle{alpha}
\bibliography{2DFPUNLSbib}


\end{document}